\RequirePackage{fix-cm}
\documentclass{svjour3}                     
\smartqed  
\newcommand{\be}{\begin{equation}} 
\newcommand{\ee}{\end{equation}}
\newcommand{\beq}{\begin{eqnarray}}
\newcommand{\eeq}{\end{eqnarray}}
\newcommand{\nbeq}{\begin{eqnarray*}}
\newcommand{\neeq}{\end{eqnarray*}}

 \journalname{Methodology and Computing in Applied Probability}

\begin{document}

\title{Time to extinction in subcritical two-sex branching processes
}


\author{David M. Hull        \and
        Manuel Mota \and George P. Yanev
}


\institute{  David M. Hull  \at
              Department of Mathematics and Computer Science, Valparaiso University, 
              Valparaiso, IN, USA \\
              \email{david.hull@valpo.edu}
           \and
           Manuel Mota \at
Department of Mathematics, University of Extremadura, Badajoz, Spain\\
              \email{mota@unex.es}
              \and
              George P. Yanev \at
              Department of Mathematics, The University of Texas - Pan American, Edinburg, TX, USA \\
              \email{yanevgp@utpa.edu}
}
\date{Received: date / Accepted: date}

\maketitle

\begin{abstract}
Lower and upper bounds for the cumulative distribution function (cdf) of the time to extinction in a subcritical
two-sex branching process  are derived. A recursive procedure for approximating this cdf is also utilized. The results are illustrated with some simulations.
\end{abstract}

\section{Introduction}
\label{intro}

There exists a significant literature regarding extinction probabilities in bisexual Galton-Watson branching processes (BGWP).
 A detailed introduction to the BGWP model and the physical properties associated with its variables can be found in Hull (2003) and Molina (2010). Necessary and sufficient conditions for certain extinction have been known for over 25 years. This raises the question - ``If extinction is certain, when will it occur?". Here we make the first attempt to give a meaningful response to this question. Daley et al. (1986), Lemma 4.3, introduced the concept of ``stopping time" in the context of two-sex processes. A version of stopping time is used in the present work. Our approach will be to apply known results concerning time-to-extinction in the standard Galton-Watson branching process to estimate time-to-extinction probabilities in the two-sex process. Our intent is to produce the first paper on time-to-extinction for the two-sex model and thus lay a foundation for future research efforts on this subject.

Let $\{(f_{n,j},m_{n,j}): n\ge 0;\ j\ge 1\}$ be a sequence of integer-valued i.i.d. bivariate random variables. A BGWP $\{Z_n\}_{n\ge 0}$ is defined by the recursion:
\beq \label{def1}
Z_0 & = & i\ge 1, \nonumber\\
Z_{n} & = & \zeta\left( \sum_{j=1}^{Z_{n-1}}f_{n,j},\ \sum_{j=1}^{Z_{n-1}}m_{n,j}\right), \quad n\ge 1,
\eeq
where $\zeta(x, y)$ is a deterministic {\it mating} function. We assume
 \begin{itemize}
      \item [(A1)] $\zeta(x,y)$ is superadditive, i.e., for any $x_1$, $x_2$, $y_1$ and $y_2\in [0,\infty)$
          $$
          \zeta(x_1+x_2,y_1+y_2)\ge \zeta(x_1,y_1)+\zeta(x_2,y_2),
          $$
      \item [(A2)] $\zeta(x,y)\le x$.
    \end{itemize}
 Both hypotheses are natural, at least from a population dynamics outlook. For an intuitive interpretation of assumption (A1), think of a two-sex population where all of the males and females are able to communicate and interact with one another without class distinctions. Superadditivity implies that the number of mating units in that scenario will not be smaller than in a scenario where the population is partitioned into a number of non-communicating groups and the mating takes place in each of these groups separately. The assumption (A2) reflects the fact that in many human and animal populations a female is allowed only one mate while a male may mate with several females. While this is not universally true, it is a common practice due to male dominance and the greater effort that females generally must make in the reproduction process. Notice also that the equality $\zeta(x,y)=x$ in assumption (A2) yields the standard asexual Galton-Watson process.

  Daley (1968)  suggested two mating functions relevant to human and animal mating. He called the first {\it completely promiscuous mating},  $\zeta(x, y) = x\min\{1,y\}$, where $x$ is the number of females and $y$ is the number of males in a given generation. In this case a {\it champion male} arises in each generation and has the capabilities to mate with every female in that generation. All other males are not allowed to mate. The other mating function is {\it polygamous mating with perfect fidelity},  $\zeta(x, y) =\min\{x, ky\}$, where $k$ is a positive integer. It is assumed that each individual will mate if a mate is available. Females may have no more than one mate and  males may have up to $k$ mates. Both mating functions, as well as most mating functions considered in the literature, satisfy (A1) and (A2).

Let $T:=\min\{n\ge 1: Z_n=0\}$ be the {\it time-to-extinction} in a BGWP. That is, $T$ counts the number of generations up to and including the generation in which extinction (absorption at zero) occurs.
Let us call $T_a$ the time to extinction in the standard (asexual) Galton-Watson process with $i\ge 1$ ancestors. Due to the independence of lines of descent (additive property, Athreya and Kaplan 1978), we have $T_a = \max_{1\le j\le i} T_a^{(j)}$, where $T_a^{(j)}$ is the time to extinction of the line of descent of the $j$th ancestor. Denote $P_i(\cdot):= P(\cdot \ |\ i\ \mbox{ancestors})$. Since $T_a^{(j)}$, $1\le j\le i$, are independent and identically distributed, we have
     $P_i(T_a \le n) = [P_1(T_a \le n)]^i$. The presence of mating in the BGWP complicates the analysis. In particular, the additive property does not hold true. As a result some properties of $T_a$ do not carry over to $T$. For instance, as a direct consequence of Theorem~3 in Gonzalez and Molina (1997), we have
     $P_i(T \le n) \le [P_1(T \le n)]^i$ for $i\ge 1$  and $n\ge 0$.

  Define the mean growth rates $r_i:=i^{-1}E[Z_{n+1}|Z_n=i]$ for $i\ge 1$. For an asexual process with offspring mean $\mu$, $r_i=\mu$ for all $i$. It is well-known that $P_i(T_a<\infty)=1$ if and only if $\mu\le 1$. In contrast, a BGWP does not have a clear genealogical tree and the changes in the population size depend on both reproduction and formation of couples (mating units). Daley et al. (1986) showed that, assuming superadditivity,  $P_i(T<\infty)=1$ for any $i\ge 1$ if and only if
$r:=\sup_{i\ge 1}r_i\le 1$. Hence, we call a BGWP subcritical, critical, and supercritical when $r<1$, $r=1$, and $r>1$, respectively.

    Our goal is to find upper and lower bounds for  $P_i(T\le n)$ and $E_i[T]:=E[T|Z_0=i]$ in the subcritical case $(r<1)$.
    In Section~2 we utilize two standard processes for estimating $P_i(T\le n)$ and derive inequalities for the expected value of $T$. Section~3 presents an alternative recursive algorithm for calculating upper and lower bounds of $P_i(T\le n)$.  In Section~4, our results are illustrated with one example and simulations.  Finally, the last section summarizes the paper's findings and lists some topics for further consideration.

\section{Asexual process bounds}
\label{sec:1}

 Define a {\it females generate females} asexual process (FGFP) $\{F_n\}_{n\ge 0}$  by:
\beq \label{def2}
F_0 & = & i\ge 1 \nonumber \\
F_{n} & = & \sum_{j=1}^{F_{n-1}} f_{n,j},\quad n\ge 1.
\eeq
Assuming (A2), it is not difficult to prove by induction that for every $n\ge 0$
\be \label{stoch1}
Z_{n}\le F_n.
\ee
 Indeed, $Z_0=i=F_0$ by definition and  for $n\ge 1$, we have
\[
Z_{n}=\zeta\left(\sum_{j=1}^{Z_{n-1}} f_{n,j}, \sum_{j=1}^{Z_{n-1}} m_{n,j}\right) \le \sum_{j=1}^{Z_{n-1}} f_{n,j} \le \sum_{j=1}^{F_{n-1}} f_{n,j}=F_n,
\]
where the last inequality follows by the induction hypothesis.

{\bf Remark.}\ For a more intuitive argument supporting (\ref{stoch1}), note that in any realization of BGWP (\ref{def1}) and FGFP (\ref{def2}), the number of reproducing females in the latter will not be smaller than the number of reproducing mating units in the former.

If $\mu_f:=E[f_{1,1}]$, then inequality (A2) implies
\[
\mu_f =  \frac{1}{i}E\left[\sum_{j=1}^i f_{n+1,j}\right]
 \ge  \frac{1}{i}E\left[\zeta\left(\sum_{j=1}^if_{n+1,j},\sum_{j=1}^im_{n+1,j}\right)\right]
 =  r_i.
\]
 Therefore,
\be \label{moment_ineq1}
r=\sup_{i\ge 1}r_i\le \mu_f.
\ee
Define a second asexual process $\{S_n\}_{n\ge 1}$ by the recursion:
\beq \label{def12}
S_0 & = & i\ge 1, \nonumber \\
S_{n} & = & \sum_{j=1}^{S_{n-1}} h_{n,j},\quad n\ge 1,
\eeq
 with $h_{n,j}=\zeta(f_{n,j},m_{n,j})$. Assuming superadditivity of $\zeta(x,y)$, one can prove that for all $n\ge 0$
 \be \label{stoch2}
 Z_n\ge S_n .
 \ee
Indeed, $S_0=i=Z_0$ and assuming $Z_{n-1}\ge S_{n-1}$ for $n\ge 1$, we have
\[
Z_n  =   \zeta\left(\sum_{j=1}^{Z_{n-1}}f_{n,j},\sum_{j=1}^{Z_{n-1}}m_{n,j}\right)\ge \sum_{j=1}^{Z_{n-1}}\zeta\left(f_{n,j},m_{n,j}\right)\ge
\sum_{j=1}^{S_{n-1}}\zeta\left(f_{n,j},m_{n,j}\right)=S_n.
\]

{\bf Remark.}\ The process (\ref{def12}) can be also viewed as a BGWP with mating function $\zeta(x,y)$ subject to the additional restriction that a male and a female may mate only if they were generated by the same mating unit, i.e., are siblings. We call this {\it siblings mating only process} (SMOP). The superadditivity assumption ensures that there will not be an advantage (i.e., more mating units) if sibling-only mating is required. This is reflected in inequality (\ref{stoch2}).

Denote
 $ \mu_s:=E[h_{1,1}]=E[\zeta(f_{n,j},m_{n,j})]$.
If $\zeta$ is superadditive, then
\[
r_i  =   \frac{1}{i}E\left[\zeta\left( \sum_{j=1}^if_{n,j},\ \sum_{j=1}^im_{n,j}\right)\right]
     \ge  \frac{1}{i} \sum_{j=1}^i E[\zeta(f_{n,j},m_{n,j})]
     =  \mu_s.
    \]
     Therefore,
\be \label{moment_ineq2}
r=\sup_{i\ge 1}r_i\ge \mu_s.
\ee

Given a BGWP, call the asexual processes $\{F_n\}$ and $\{S_n\}$  its ``associated processes".

Next we recall known (Agresti 1974, Theorem~1 and Haccou et al. 2005, Theorem~5.3) inequalities for the cdf of the time-to-extinction in an asexual process. For any standard process with offspring law $\{p_k\}_{k\ge 0}$, offspring mean $0<\mu<1$ and offspring variance $\sigma^2<\infty$ define
\[
\overline{Q}_i(\mu, p_0):= 1-\frac{1-\mu^n}{c_1(1-\mu)+1-\mu^n}\delta_{i,1},
\]
where $c_1=\max\left\{2, \  \mu(\mu+p_0-1)^{-1}\right\}$ and $\delta_{1,1}=1$; $\delta_{i,1}=0$ for $i\ne 1$.
Define also
\[
\underline{Q}_i(\mu, \sigma^2):=\frac{c_2}{1-(1-c_2)\mu^n}-\frac{1-\mu^n}{1-(1-c_2)\mu^n}\frac{i-1}{2},
\]
where  $c_2=\mu(1-\mu)\sigma^{-2}$. If $\mu<1$ and $\sigma^2<\infty$, then
\be \label{ineq}
\underline{Q}_i(\mu, \sigma^2)\mu^n \le P_i(T_a>n)\le \overline{Q}_i(\mu, p_{0})i\mu^n,\qquad n\ge 0.
\ee
Denote $p_{f,0}:=P(f_{1,1}=0)$.

\vspace{0.3cm}
{\bf Proposition 1}\   Assume (A1) and (A2). If $\mu_f:=E[f_{1,1}]<~1$ and $\sigma_s^2:=Var[\zeta(f_{1,1},m_{1,1})]<\infty$, then for $i
\ge 1$
\be \label{ineq_main}
\underline{Q}_i(\mu_s, \sigma_s^2)\mu_s^n \le P_i(T>n)\le \overline{Q}_i(\mu_f, p_{f,0})i\mu_f^n,\qquad n\ge 0,
\ee
and
\be \label{distance}
\Delta_i(n):=\overline{Q}_i(\mu_f, p_{f,0})i\mu_f^n-\underline{Q}_i(\mu_s, \sigma_s^2)\mu_s^n =O\left(\mu^n_f\right), \qquad n\to \infty.
\ee
\noindent{\bf Proof}.
  Let $T_{FGFP}$ and $T_{SMOP}$ denote the times to extinction in FGFP
  and SMOP, respectively. The inequalities (\ref{stoch1}) and (\ref{stoch2}) imply for all $n\ge 0$
\be \label{processes_ineq1}
P_i(T_{SMOP}> n)\le P_i(T>n)\le P_i(T_{FGFP}> n).
\ee
Applying (\ref{ineq}) to $P_i(T_{SMOP}>n)$ and $P_i(T_{FGFP}> n)$, we obtain
(\ref{ineq_main}).

Since  $\mu_s\le \mu_f$ and
\[
\Delta_i(n)=\mu_f^n\left[\overline{Q}_i(\mu_f, p_{f,0})i-\underline{Q}_i(\mu_s, \sigma_s^2)\left(\frac{\mu_s}{\mu_f}\right)^n\right],
\]
taking into account that
\[
\lim_{n\to \infty}\overline{Q}_i(\mu_f, p_{f,0})=1-\frac{\delta_{i,1}}{c_1(1-\mu_f)+1}\quad \mbox{and}\quad \lim_{n\to \infty}\underline{Q}_i(\mu_s, \sigma_s^2)=c_2-\frac{i-1}{2},
\]
where $c_1=\max\left\{2, \  \mu_f(\mu_f+p_{f,0}-1)^{-1}\right\}$  and $c_2=\mu(1-\mu_s)\sigma_s^{-2}$, we obtain (\ref{distance}). The proof is complete.

{\bf Remarks.}\ The assumptions in Proposition~1 are common in the literature on two-sex branching processes. We have already discussed (A1) and (A2) in Section~1. Furthermore, small numbers of mating units and hence eventual extinction is primarily caused by small numbers of females in the various generations. Small numbers of males can also cause extinction  but that can be dealt with by allowing males to have multiple mates. Hence, the assumption that $E[f_{1,1}] < 1$. The assumption $\sigma_s^2 < \infty$ is needed to apply Agresti's bound. Note that for specific mating functions the last two hypotheses simplify as follows. If $\zeta(x,y)=\min\{x,ky\}$ then $r=\min\{E[f_{1,1}],kE[m_{1,1}]\}$ and therefore $E[f_{1,1}]<1$ implies $r<1$. If $\zeta(x,y)=x\min\{1,y\}$ then $r=E[f_{1,1}]$.
Also, since $\zeta(x,y)\le x$, $E[\zeta(f_{1,1},m_{1,1})^2]\le E[f_{1,1}^2]$ and therefore $Var[f_{1,1}]<\infty$ implies $\sigma_s^2<\infty$.

Next we will find lower and upper bounds for the expected time to extinction, $E_i[T]$. For the subcritical $(\mu<1)$ asexual process (Haccou  et al. 2005, Thm.~5.4)
\be \label{a_ineq}
\left(\frac{\ln i-\ln \ln i}{|\ln \mu |}-1\right)\left(1-\frac{1}{ib_a}\right)\le E_i[T_a]\le \frac{\ln i}{|\ln \mu|}+\frac{2-\mu}{1-\mu},
\ee
where $i\ge 3$ and $b_a\ge (1-\mu)\mu/\sigma^2$, provided that the offspring variance $\sigma^2<\infty$.  The following result holds true for a BGWP.

\vspace{0.3cm}
{\bf Proposition 2}\ Assume (A1), (A2), and $\sigma_s^2:=Var[\zeta(f_{1,1},m_{1,1})]<~\infty$. If $0<r<1$, then for $i\ge 3$
\be  \label{mean_time}
\left(\frac{\ln i-\ln \ln i}{|\ln \mu_s |}-1\right)\left(1-\frac{1}{ib_s}\right)\le E_i[T]\le \min\left\{\frac{\ln i}{|\ln r|}+\frac{2-r}{1-r},\ \frac{i}{1-r}\right\},
\ee
where $b_s\ge (1-\mu_s)\mu_s/\sigma^2_s$.

{\bf Proof.}
Since
$
P_i(T>n)=P_i(Z_n\ge 1)\le E_i[Z_n]\le ir^n$,
we have
\be \label{upper1}
E_i[T] = \sum_{n=0}^\infty P_i(T>n)\le i\sum_{n=0}^\infty r^n=\frac{i}{1-r}.
\ee
On the other hand,
\beq \label{upper}
E_i[T] & = & \sum_{n=0}^\infty P_i(T>n) \\
    & \le &
\sum_{0\le n< \ln i/|\ln r|}P_i(T>n)+i\sum_{n\ge \ln i/|\ln r|}r^n \nonumber \\
    & \le &
\frac{\ln i}{|\ln r|}+1+\frac{ir^{\ln i/|\ln r|}}{1-r} \nonumber \\
    & = & \frac{\ln i}{|\ln r|}+\frac{2-r}{1-r}. \nonumber
\eeq
Equations (\ref{upper1}) and (\ref{upper}) imply the upper bound in (\ref{mean_time}).
Let us now find a lower bound for $E_i[T]$. Recalling the SMOP and referring  to (\ref{moment_ineq2}), (\ref{processes_ineq1}), and (\ref{upper1}) we obtain
\be \label{lower}
E_i[T] =  \sum_{n=1}^\infty P_i(T>n)
     \ge
 \sum_{n=1}^\infty P_i(T_{SMOP}>n)
     =
E_i[T_{SMOP}],
\ee
where $T_{SMOP}$ denote the time to extinction of SMOP (\ref{def12}). Equation (\ref{lower}), in view of (\ref{a_ineq}), implies the left-hand side inequality in (\ref{mean_time}), which completes the proof.

\vspace{0.3cm}
{\bf Analytical Example}\ Write $\xi\in MG(b,c)$ when the random variable $\xi$ follows the modified geometric distribution with parameters $b>0$ and $c>0$ such that $b+c\le 1$ and probability mass function
\[
P(\xi=k)=bc^{k-1}\quad \mbox{for}\quad k=1,2,\ldots; \quad P(\xi=0)=1-\sum_{k=1}^\infty P(\xi=k).
\]
Define a BGWP $\{Z_n\}$ with $Z_0=1$ and for $i,j\ge 0$
\[
P(f_{1,1}=i, m_{1,1}=j)=P(f_{1,1}=i)P(m_{1,1}=j),
\]
where $f_{1,1}\in MG(b_f, c_f)$ and $m_{1,1,}\in MG(b_m,c_m)$, i.e., the numbers of both female and male offspring follow independent modified geometric distributions.
Assume also $\zeta(x,y)=x\min \{1,y\}$. For the time to extinction of the associated FGFP we have (e.g., Taylor and Karlin (1994), p.457)
\[
P(T_{FGFP}> n)=\left(\frac{u_f-1}{u_f-\mu_f^n}\right)\mu_f^n,
\]
where $\mu_f=b_f/(1-c_f)^2$ and $u_f:=(1-b_f-c_f)/(c_f(1-c_f))>1$ if $\mu_f<1$. On the other hand, for the associated SMOP,
$
P(h_{1,1}=0)=1-P(f_{1,1}>0)P(m_{1,1}>0)$ and $P(h_{1,1}=k)=P(f_{1,1}=k)P(m_{1,1}>0)$
for $k=1,2,\ldots$ It is not difficult to verify that $h_{1,1}\in MG(b_s, c_f)$, where $b_s=b_fb_m/(1-c_m)$. Therefore,
\[
P(T_{SOMP}>n)=\left(\frac{u_s-1}{u_s-\mu_s^n}\right)\mu_s^n=\left(\frac{u_s-1}{u_s-\mu_s^n}\right)\left(\frac{b_m}{1-c_m}\right)^n \mu_f^n,
\]
where $u_s:=(1-b_s-c_f)/(c_f(1-c_f))>1$ if $\mu_s<1$.
Assuming $\mu_f<1$, the inequalities (\ref{processes_ineq1}) become
\be \label{ex_ineq}
\mu_f^n\left(\frac{b_m}{1-c_m}\right)^n\left(\frac{u_s-1}{u_s-\mu_s^n}\right)\le P_1(T>n)\le \mu_f^n\left(\frac{u_f-1}{u_f-\mu_f^n}\right)
\ee
Note that the maximum error of estimating $P_1(T>n)$ is $O(\mu_f^n)$ as $n\to \infty$. Finally, using (\ref{ex_ineq}) we obtain for $\mu_f<1$
\[
\left(1-\frac{1}{u_s}\right)\frac{1}{1-\mu_s}\le E_1[T]\le \frac{1}{1-\mu_f}.
\]

{\bf Remark.}\ The procedure for estimating $P_i(T\le n)$ presented in this section can be extended to the cases when the associated process $\{F_n\}$ (and also $\{S_n\}$) is not subcritical, i.e., $\mu_f\ge 1$.
If $\{F_n\}$ is critical, then (\ref{ineq}) can be replaced by the interval formulas in Theorem~1 of Agresti (1974). When $\{F_n\}$ (and maybe also $\{S_n\}$) is supercritical, approximating intervals for the distribution of time-to-extinction can be derived by using the duality between subcritical and supercritical asexual processes. The theory behind this duality is discussed in Agresti (1974), see also Athreya and Ney (1972), pp. 72-73. In particular, it is well-known that if $T_a^{sup}$ and $T_a^{sub}$ are the times to extinction in a supercritical and its dual subcritical asexual process, respectively, then $a < P(T_a^{sub} < n) < b$ implies $qa < P(T_a^{sup} < n) < qb$, where $q$ is the extinction probability of the supercritical process with one individual in the initial generation

\section{Finite Markov chain approximations}
\label{sec:2}

 The inequalities (\ref{ineq_main}) are intended to be the basis of estimating a time-to-extinction probability. If this approach lacks precision, i.e., the interval lengths for various values of $n$ are too wide, then one can use a recursive algorithm for approximating $P_i(T\le n)$, which will be developed in this section. The procedure is based on two finite Markov chains. For a fixed integer $M\ge i$, let $n(M):=\min \{n: Z_n>M\}$. For $n\ge 1$ define

\vspace{0.5cm}$
\widetilde{Z}_n=
\left\{
  \begin{array}{ll}
    Z_n & \mbox{if}\quad n<n(M), \\
    0 & \mbox{if}\quad n\ge n(M),
  \end{array}
\right.
$ \hspace{0.5cm}and \hspace{0.5cm}
$
\widehat{Z}_n=
\left\{
  \begin{array}{ll}
    Z_n & \mbox{if}\quad n<n(M), \\
    M & \mbox{if}\quad n\ge n(M).
  \end{array}
\right.
$

\vspace{0.5cm} \noindent
Note that $\{\widehat{Z}_n\}$ appeared in Daley et al. (1986).

Let $\widetilde{T}$ and $\widehat{T}$ be the times to extinction for $\{\widetilde{Z}_n\}_{n\ge 0}$ and $\{\widehat{Z}_n\}_{n\ge 0}$, respectively.   Define  $\widetilde{G}_{i}(n):=P_i(\widetilde{T}\le n)$, $\widehat{G}_{i}(n):=P_i(\widehat{T}\le n)$, and $G_{i}(n):=P_i(T\le n)$.
By the definitions of $\{\widetilde{Z}_n\}$ and $\{\widehat{Z}_n\}$, we have for $n\ge 0$
\[
P_i(\widehat{Z}_n=0)\le P_i(Z_n=0)\le P_i(\widetilde{Z}_n=0),
\]
which is equivalent to
\be \label{tilde_hat}
\widehat{G}_{i}(n)\le G_i(n)\le \widetilde{G}_{i}(n).
\ee
Define $P_{ij}:=P(Z_{n+1}=j|Z_n=i)$ for nonnegative integers $i, j$, and $n$. The utility of (\ref{tilde_hat}) is that it is not difficult to compute the bounds $\widehat{G}_{i}(n)$ and $\widetilde{G}_{i}(n)$ for any $n\ge 0$ using the recurrences: $\widehat{G}_{i}(0)=0$,
\be \label{rec1}
\widehat{G}_{i}(n)=P_{i0}+\sum_{j=1}^M P_{ij}\widehat{G}_{j}(n-1)
\ee
and $\widetilde{G}_{i}(0)=0$,
\beq  \label{rec2}
\widetilde{G}_{i}(n) & = & P_{i0}+P_i(Z_1>M)+\sum_{j=1}^M P_{ij}\widetilde{G}_{j}(n-1)\\
    & = & 1 - \sum_{j=1}^M P_{ij}(1-\widetilde{G}_{j}(n-1)). \nonumber
\eeq
The next result 
shows, in particular, that $\widetilde{G}_{i}(n)- \widehat{G}_{i}(n)=O(1/M)$ as $M\to \infty$.

\vspace{0.3cm}
{\bf Proposition 3} If $0<r< 1$, then for $i\ge 1$ and $n\ge 2$,
\be \label{eq21}
0\le \widetilde{G}_{i}(n)-G_{i}(n)\le \frac{ic_{n+1}(r)}{M}\quad \mbox{and}\quad  0\le G_{i}(n)-\widehat{G}_{i}(n)\le
\frac{ic_n(r)}{M},
\ee
 where $c_n(r)=(r-r^n)/(1-r)$.

\noindent{\bf Proof.}\
For the first inequality in (\ref{eq21}), observe that
\nbeq
P_i(\widetilde{Z}_n=0) & = & P_i\left(\bigcup_{k=1}^n \{Z_k=0\}\right)+P_i\left(\bigcup_{k=1}^n \{Z_k>M\}\right)\\
    & \le & P_i(Z_n=0)+\sum_{k=1}^n P_i(Z_k>M).
\neeq
Therefore,
\be \label{first_inequality}
\widetilde{G}_{i}(n)-G_{i}(n)\le \sum_{k=1}^{n}P_i(Z_k>M).
\ee
In addition, since $E[Z_k |Z_{k-1}]=Z_{k-1}r_{Z_{k-1}}$ for a positive $k$, we have $E[Z_k]= E[E[Z_k|Z_{k-1}]]= E[Z_{k-1}r_{Z_{k-1}}]\le  r E[Z_{k-1}]$, which after iterating implies
\be\label{expectation_inequality}
E_i[Z_k]
     \le
        ir^k.
\ee
The inequalities (\ref{first_inequality}) and (\ref{expectation_inequality}), along with the Markov inequality, yield
\be \label{Markov_ineq}
\widetilde{G}_{i}(n)-G_{i}(n)\le \sum_{k=1}^{n}P_i(Z_k>M)
     \le  \frac{1}{M}\sum_{k=1}^{n}E_i[Z_k]
     \le
\frac{i}{M}\sum_{k=1}^{n}r^k.
\ee
From (\ref{tilde_hat}) and (\ref{Markov_ineq}) we obtain the first part of (\ref{eq21}). To prove the second inequality in (\ref{eq21}), we will first show by induction that
\be \label{eq2a}
G_{i}(n)-\widehat{G}_{i}(n)\le
\sum_{k=1}^{n-1}P_i(Z_k>M).
\ee
The following recursion holds true
\beq \label{rec_M}
G_i(0)& = & 0, \nonumber \\
G_{i}(n)
     & = &   P_{i0}+\sum_{j=1}^{\infty}P_{ij}G_{j}(n-1), \quad n\ge 1.
\eeq
From (\ref{rec_M}), we have for $n=1$ and  $M\ge i$,
\nbeq
G_{i}(2)
    & = &
        P_{i0}+\sum_{j=1}^M P_{ij}G_{j}(1)+\sum_{j=M+1}^\infty P_{ij}G_{j}(1)\\
        & \le &
        \widehat{G}_{i}(2)+P_i(Z_1>M),
\neeq
which establishes (\ref{eq2a}) for $n=2$.  Assume (\ref{eq2a}) is true for a fixed $n>2$. Then for any $M\ge i$,
\nbeq
\lefteqn{G_{i}(n+1)
     =
P_{i0}+\sum_{j=1}^M P_{ij}G_{j}(n)+\sum_{j=M+1}^\infty P_{ij}G_{j}(n)}\\
    & \le &
P_{i0}+\sum_{j=1}^MP_{ij}\left[ \widehat{G}_{j}(n)+\sum_{k=1}^{n-1}P_j(Z_k>M)\right]+P_i(Z_1>M)\\
    & = &
    \widehat{G}_{i}(n+1)+\sum_{j=1}^MP_{ij}\left[\sum_{k=1}^{n-1}P_j(Z_k>M)\right]+P_i(Z_1>M)\\
    & \le &
    \widehat{G}_{i}(n+1)+\sum_{k=1}^{n-1}\sum_{j=1}^\infty P_{ij}P_j(Z_k>M)+P_i(Z_1>M)\\
    & = &
    \widehat{G}_{i}(n+1)+\sum_{k=1}^{n-1}P_i(Z_{k+1}>M)+P_i(Z_1>M)\\
    & = &
    \widehat{G}_{i}(n+1)+\sum_{k=1}^{n}P_i(Z_{k}>M),
    \neeq
which completes the proof of (\ref{eq2a}). The second part of (\ref{eq21}) follows from (\ref{tilde_hat}) and (\ref{eq2a}) applying the Markov inequality as in (\ref{Markov_ineq}). The proposition is proved.

\vspace{-0.3cm}\section{Numerical Example}
\label{sec:4}

In this section we illustrate through an example the proposed approximation procedures.
Consider a BGWP with $Z_0 = i\ge 1$ and for
 $0<\alpha <1$
\[
P(f_{1,1}=k, m_{1,1}=j)=\left\{
                          \begin{array}{ll}
                           {3 \choose k}\alpha^k (1-\alpha)^j & \mbox{for}\ k+j=3, \ k,j\ge 0,\\
                            \ 0 & \mbox{otherwise}.
                          \end{array}
                        \right.
\]
Assume also a completely promiscuous mating function $\zeta (x,y)=x\min\{1,y\}$.
The distribution of the number of females generated by a mating unit is
$p^f_{k}:={3 \choose k}\alpha^k(1-\alpha)^{3-k}$, $0\le k\le 3$. It follows from (\ref{moment_ineq1}) that the process is subcritical provided $E[f_{1,1}]=3\alpha<1$. We set $\alpha=0.25$.

To obtain asexual process bounds for $P_i(T\le n)$, let us consider the associated processes  FGFP and SMOP defined by (\ref{def2}) and (\ref{def12}), respectively, with offspring pgf $g(x)=(\alpha x+1-\alpha)^3$ and $h(x)=3\alpha^2(1-\alpha)x^2+3\alpha(1-\alpha)^2x+\alpha^3+(1-\alpha)^3$.  If $f_n(x)$ denotes the $n$th composition of a function $f(x)$ with itself, then it is well-known that
\be \label{a_pgfs}
P_i(T_{FGFP}\le n)=(g_n(0))^i \quad \mbox{and}\quad   P_i(T_{SMOP}\le n)=(h_n(0))^i.
\ee
Table 1 lists $G^{FGFP}_i(n):=P_i(T_{FGFP}\le n)$ and $G^{SMOP}_i(n):=P_i(T_{SMOP}\le n)$ for $\alpha=0.25$ and selected values of $i$ and $n$. Here we use the explicit formulas (\ref{a_pgfs}) instead of the inequalities (\ref{ineq}).


Next, we compute $\hat G_i(n)$ and $\tilde G_i(n)$ using the recurrences (\ref{rec1}) and (\ref{rec2}) with $M=20$. The one-step transition probabilities $P_{ij}, j=0,\dots M$ have been estimated by the Monte Carlo method based on 10,000 simulations of $Z_1$ when $Z_0=i$. The obtained bounds for $P_i(T\le n)$ agree up to the third decimal place for all $i$ and $n$ in Table~1, where this common value is denoted as $G_i^*(n)$.

The results of both estimating procedures described above are presented in Table~1.

\begin{table}
 \caption{Comparison between $G^{FGFP}_i(n)$, $G^*_i(n)$, and $G^{SMOP}_i(n)$}
\label{tab:1}
 \begin{tabular}{lllll}
\hline\noalign{\smallskip}
  $i$  & $n$  & $G^{FGFP}_i(n)$ & $G^*_i(n)$ & $G^{SMOP}_i(n)$ \\
\noalign{\smallskip}\hline\noalign{\smallskip}
 2 & 2& 0.392 & 0.400 & 0.421 \\
 \noalign{\smallskip}\hline\noalign{\smallskip}
 & 5 &  0.759 & 0.773 & 0.801  \\
   \noalign{\smallskip}\hline\noalign{\smallskip}
 & 7 &  0.868 &  0.880 &  0.902 \\
 \noalign{\smallskip}\hline\noalign{\smallskip}
& 10&  0.945 &  0.952 &   0.966 \\
\noalign{\smallskip}\hline \hline\noalign{\smallskip}
5  & 2 & 0.096  & 0.098  &  0.115  \\
\noalign{\smallskip}\hline\noalign{\smallskip}
  & 5 & 0.503 & 0.513 &  0.575 \\
  \noalign{\smallskip}\hline\noalign{\smallskip}
  & 7 &  0.702 &  0.714 &  0.774 \\
  \noalign{\smallskip}\hline\noalign{\smallskip}
  &  10 &  0.869 &  0.879 &   0.912 \\
  \noalign{\smallskip}\hline \hline \noalign{\smallskip}
 10 & 2 & 0.009&  0.009 &  0.013  \\
 \noalign{\smallskip}\hline\noalign{\smallskip}
  &  5 & 0.253  & 0.259  &  0.330  \\
  \noalign{\smallskip}\hline\noalign{\smallskip}
   & 7 & 0.492 & 0.503  &   0.598 \\
   \noalign{\smallskip}\hline\noalign{\smallskip}
  &  10 &  0.755 &  0.767 &  0.842 \\
  \noalign{\smallskip}\hline
\end{tabular}
 \end{table}

\vspace{-0.3cm}\section{Concluding remarks}
\label{sec:5}

It is our hope and expectation that this paper will provide a starting point and motivation for further consideration of the time to extinction in both theoretical context and applications to certain traits or characteristics of two-sex species. Such may include, for instance, various genotypes of intermediate forms in the evolution of a particular two-sex species.

To summarize, two kinds of bounds for the distribution function of $T$ are derived. The first is based on constructing ordinary branching processes which stochastically bound the process specified by the number of mating units in each generation. The (quite mild) moment assumptions give bounds for $P_i(T>n)$ of order $\mu^n$, where $\mu$ is the per-capita mean number of offspring of the bounding processes. These means differ, but the bounds show that the upper bounding mean is $\mu_f$, the mean number of female offspring in the two-sex process, assumed less than unity. The bounds are specified in terms  of offspring means and variances, hence easy to compute and also $P_i(T>n)=O(\mu_f^n)$ as $n\to \infty$.

The second pair of bounds is based on two sequences of finite-state Markov chains obtained by stopping the two-sex branching process in two different ways. The zero hitting times of these Markov chains bound the extinction time of the two-sex process. Each sequence is indexed by a positive integer $M$, such that the corresponding chain has state-space $\{0, 1, \ldots, M\}$.
An estimate of the error involved with each bound is computed. It has a very simple form and is proportional to $M^{-1}$. Thus, $M$ can be chosen to achieved a prescribed maximum error. These bounds require numerical values of the one-step transition probabilities and much more computation that the first ones. However, they can be made as tight as desired.

It would be more efficient to deal with probabilities involving $T$ when the mating function is further specified. For instance, the counter-example given in Hull (1982) does not require estimates and $P(T=n)$ can be calculated directly. Are there other mating functions with this property? Also, Daley's relevant mating functions should be considered as separate projects when investigating probabilities of $T$. Are there shortcuts when dealing with these specific mating functions?

\begin{acknowledgements}
The authors thank the anonymous referee for the valuable comments and suggestions.
\end{acknowledgements}

\end{document}